\documentstyle{amsppt}
\magnification=1100

\input xy
\xyoption{all}

\topmatter
\title Explicit Constructions for Genus 3 Jacobians \endtitle
\author Jesus Romero-Valencia \& Alexis G. Zamora \endauthor
\address Universidad de Guerrero \newline Av. L\'azaro C\'ardenas
s/n
\newline Ciudad Universitaria, C.P. 39087\newline Chilpancingo,
Guerrero, M\'exico
\endaddress \email  jromv\@ yahoo.com\endemail
\address U. A. Matem\'aticas, U.de  Zacatecas\newline
Camino a la Bufa y Calzada Solidaridad, C.P. 98000\newline
Zacatecas, Zac. M\'exico\endaddress \email alexis\@ mate.reduaz.mx
\endemail

\abstract Given a canonical genus three curve $X=\{ F=0\}$, we
construct, emulating Mumford discussion for hyperelliptic curves,
a set of equations for an affine open subset of the jacobian $JX$.
 We give explicit algorithms describing the law
group in $JX$. Finally we introduce a related construction by
means of an imbedding of the open set previously described in a
Grassmanian variety.
\endabstract
\thanks The authors were partially supported by CONACyT Grant
25811\endthanks \subjclass 14H40, 14H45 \endsubjclass
\endtopmatter

\document

\head 1. Introduction \endhead

The aim of this paper is to construct, in the spirit of [7], a set
of explicit equations for an affine open set of the Jacobian
Variety of a canonical genus three curve $X$, defined by a quartic
equation $F=0$. We develop in details the case when $X$ has a
hyper-flex point, but the construction is valid for any genus $3$
non-hyperelliptic curve. By a hyper-flex we understand a point
$p\in X$ such that its tangent line intersects $X$ of order $4$ at
$p$. Our ground field is $\Bbb C$, but the results seem to be
valid over any algebraically closed field of characteristic
different from $2$ or $3$.

The idea is to work with a subset of degree $3$ non-special
effective divisors and associate to it a pencil of conics. Not
only an open set of $JX$ can be described, but also the
multiplication by $-1$ and the law of group (even when this last
not explicitly in terms of affine coordinates).

In section 3 we construct an affine variety $Z\subset \Bbb P^{17}$
that corresponds bijectively with a subset of $Pic^3(X)$. $Z$ is
described by means of the possible decompositions $F=AG+BH$ for
certain class of conics $A$ and $B$. We prove that this
correspondence is in fact an isomorphism of algebraic varieties.

The fact that points in an open set of the Jacobian can be put in
correspondence with a decomposition of the type $F=AG+BH$ is by no
means surprisingly, this follows from a homological construction,
namely the computation of a minimal resolution of the $\Cal
O_{\Bbb P^2}-$module $\Cal O_X(D-3\infty)$ for $D$ a general
element of $Pic^3(X)$ (and thus $[D-3\infty]$ a general element of
$Jac X$). The degeneracy locus of this resolution is, by one hand
the curve $X$ itself, and by the other the determinant of a
$2\times 2$ matrix with coefficients in $H^0(\Bbb P^2,\Cal O_{\Bbb
P^2}(2))$. This interpretation is developed in section 2.

In section 4 we work with a curve with a hyper-flex point and we
describe multiplication by $(-1)$ in $JX$ through  a simple
geometric construction, prove that $Z$ is invariant under
multiplication by $(-1)$ and write explicitly the morphism:

$$ Z @> (-1)>> Z.$$

This description allows us to obtain equations for the open set
$Z/\{\pm 1\}$ of the Kummer variety $K(X)$. In that section we
identify $JX$ with the set of equivalence class of divisor of the
form $D-3 \infty$, with any $D$ effective divisors of degree 3 and
$\infty$ a fixed hyper-flex point of $X$.  We finish the section
by describing geometrically how to obtain the sum of two divisors
$D-3\infty$, $D'-3\infty$ $\in JX$.

Section 4 is devoted to outline a more intrinsic variant of the
previous construction. The idea is, instead to looking at the
decomposition $F=AG+BH$, to go directly to the pencil of conics
generated by $A$ and $B$ and define a morphism

$$ Z @>>> \Bbb G(1,V),$$
with $V\subset H^0(\Bbb P^2, \Cal O(2))$ the space of conics
passing through $\infty$.

Other explicit constructions of Jacobians can be found in [7],
[1], [2], [4], [9]. Undoubtedly, the beautiful construction in [7]
(that have its roots, according to Mumford, in works by Jacobi)
have been the source of inspiration for the rest of the mentioned
papers.

Even when generalizations of this construction for higher genus
could be more complicated (for instance, quotients by the action
of groups of matrices having as coefficients homogeneous
polynomials can be involved), it is clear that the homological
method of minimal resolution (see section 2) is the natural
mechanism for treat such possible generalization. The discussion
in genus $3$ is considerable simplified by the fact that the
resolution gives place to an expression of $F$ like a $2\times 2$
determinant. The authors are deeply indebted to an anonymous
referee that indicated how our original elementary construction,
developed in section 3, fits within the beautiful framework of
minimal resolutions.

\head 2 The Homological Construction \endhead

This section is based on [6]. Another important exposition of this
ideas appears in [3].

Let $\Cal F$ be a coherent sheaf on $\Bbb P^n$, then $\Cal F$ is
Cohen-Macaulay and, therefore, a minimal resolution

$$ 0 @>>> L_k @>\alpha_k >> L_{k-1} @>\alpha_{k-1} >> ... @>>> L_0 @>>> \Cal F @>>> 0,$$
with $L_p=\oplus_q (B_{p,q})\otimes\Cal O_{\Bbb P^n}(-q)$ and $k$
the codimension of $\Cal F$ in $\Bbb P^n$, exists. Here $B_{p,q}$
is a vector space whose dimension indicates the number of summands
$\Cal O_{\Bbb P^n}(-q)$ appearing in $L_p$. This minimal
resolution is unique up to endomorphism of the $L_i$ and the
$\alpha_i$'s are given by matrices of homogeneous polynomials. If
we want to compute such a resolution it is sufficient to know the
dimension of the vector spaces $B_{p,q}$.

M. Green introduced an explicit way to compute this vector spaces,
based on Koszul resolution and Euler's exact sequence. Although it
is not the only method appearing in the literature (see [3]) we
prefer here this treatment.

We are interested in a particular case: $n=2$, and $X\subset \Bbb
P^2$ is a degree $d$ curve. If $\Cal F$ is any coherent sheaf on
$X=\{ F=0 \}$ we have a minimal resolution:

$$0@>>> L_1 @>\alpha>> L_0 @>>> \Cal F@>>>0,$$
with $L_1$ and $L_2$ of the same rank. Moreover, $\det \alpha
=\lambda F$, $\lambda \in \Bbb C$, just because the vanishing of
this determinant correspond to the geometrical locus where the map
$\alpha$ is not surjective. In particular, the sum of the degrees
of the diagonal entries of $\alpha$ must equals $d$.

The spaces $B_{1,q}$ in the resolution of a coherent sheaf on $X$
can be computed in the following way:

\proclaim{Proposition 2.1} Let $\Cal O_X(D)$ be a line bundle on
$X$, then the terms $B_{1,q}$ in the minimal resolution of $\Cal
O_X(D)$ as $\Cal O_{\Bbb P^2}-$module are:

$$B_{1,q}\simeq \ker \left( H^1(X, \det \Omega_{\Bbb P^2}(1)\otimes \Cal O_X(D)(q-2)) \rightarrow \wedge^2\Bbb C^3
\otimes H^1(X, \Cal O_X(D)(q-2))\right),$$ where the maps
appearing in the statement are given by the restriction to $X$ of
the Koszul complex.

\endproclaim

\demo{Proof} This is just a combination of Theorem 2.1 and
Proposition 2.4 in [6].
\enddemo

We want to make this computation explicit in the case of a
non-singular quartic on $\Bbb P^2$, that is, a canonical genus $3$
curve. We have:

\proclaim{Proposition 2.2} Let $X=\{F=0\}$ a non-singular quartic
on $\Bbb P^2$. Let $D$ be a zero-degree divisor on $X$, $D\neq 0$,
then $\Cal O_X(D)$ admits a minimal resolution:

$$0@>>> \Cal O_{\Bbb P^2}(-3)^{\oplus 2} @>>> \Cal O_{\Bbb P^2}(-1)^{\oplus
2} @>>> \Cal O_X(D)@>>>0.$$ \endproclaim

\demo{Proof} First, we compute the dimension of the vector space
$B_{1,3}$. By the previous proposition it is the kernel of a map:

$$H^1(X, \det \Omega_{\Bbb P^2} (1)\otimes \Cal O_X(D)(1))@>>>
\wedge^2 \Bbb C^3 \otimes H^1(X,\Cal O_X(D)(1)).$$

Now, $\Cal O_X(1)\simeq \omega_X$, thus being $D$ not equivalent
to zero we have

$$h^1(X,\Cal O_X(D)(1))=h^0(X,\Cal O_X(-D))=0.$$ Thus, the dimension
of the kernel equals the dimension of the vector space

$$H^1(X,\det \Omega_{\Bbb P^2} (1)\otimes \Cal O_X(D)(1)).$$ Euler's exact
sequence implies that $\det \Omega_{\Bbb P^2}(1)=\Cal O_{\Bbb
P^2}(-1)$. Therefore,

$$H^1(X, \det \Omega_{\Bbb P^2} (1)\otimes \Cal O_X(D)(1))\simeq H^1(X,\Cal
O_X(D)),$$ and, using again that $D$ is not the zero divisor we
have that $h^1(X,\Cal O_X(D))=2$. We conclude that $\dim
B_{1,3}=2$.

Next, note that $h^1(X, \det \Omega_{\Bbb P^2} (1)\otimes \Cal
O_X(D)(k))=h^1(X,\Cal O_X(D)\otimes \omega_X^{k-1})=0$ for $k\ge
2$. Thus, no terms $B_{1,q}$ appears in the minimal resolution for
$q\ge 4$.

Finally, for terms $B_{1,q}$ with $q<3$ we note that the kernels
to be computed are those in co-homology associated to the
restricted Euler's sequence:

$$0@>>> \Omega ^1_{\Bbb P^2}(1) (D) @>>> \wedge^2 \Bbb C^3\otimes \Cal O_X(D)
@>>> \omega_X(D) @>>> 0,$$ after tensor product by
$\omega_X^{-k}$, $k>0$.

The injectivity of this maps at the 1-level co-homology follows
from the fact that $h^0(X,\omega_X^{-k}(D))=0$ for $k\le 0$.

Therefore $L_1=\Cal O_{\Bbb P^2}(-3)^{\oplus 2}$, it follows from
the fact that the sum of the degree of the diagonal entries of
$\alpha$ must be $4$ that $L_0=\Cal O_{\Bbb P^2}(-1)^{\oplus 2}$
\qed

\enddemo

In this way $\det \alpha=0$ gives rise to a decomposition
$F=AG-BH$ for some conics $(A,B,G,H)$. This expression is uniquely
determined by $D$ up to the actions by endomorphisms of $GL(2)$ on
$\Cal O_{\Bbb P^2}(-3)^{\oplus 2}$ and $\Cal O_{\Bbb
P^2}(-1)^{\oplus 2}$ respectively. In conclusion:

\proclaim{Lemma 2.2} With the same notation as before there is a
one to one correspondence between non-trivial degree-zero line
bundles on $X$ and conics $(A,B,G,H)$ satisfying $F=AG-BH$, modulo
the action of $GL(2)^{\oplus 2}$ given by:

$$ (A_0,B_0).\left( \matrix A & B \\ H & G \endmatrix \right) =
A_0^{-1} \left( \matrix A & B \\ H & G \endmatrix \right) B_0.$$
\endproclaim

The goal of the next section is to prove that this correspondence,
at least when restricted to some open set of $Jac X$, is in fact a
isomorphism of algebraic varieties, and to give an explicit
description of the quotient of quadruples of conics modulo the
action of $GL(2)^{\oplus 2}$.

Of course, we can construct a similar correspondence at least in
the case of non-singular plane curves $X$ and obtain a set of
equations that, up to the action of some algebraic group,
describe, at least birationally, an open set of $Jac X$. However,
the explicit computation of the quotient and the resulting
decomposition for the equation of $X$ could have in general not a
so simple description. For instance, if $d=5$, $g=6$ the matrix
$\alpha$ is a $5\times 5$ matrix of linear coefficients.

\head 3. Algebraic Explicit Construction for an affine set of a
genus $3$ Jacobian \endhead

Let $X=\{F=0\}$ be a plane genus $3$ curve defined by the
non-singular homogeneous quartic polynomial $F$. We assume, in
order to simplify  the discussion, that $X$ has a hyper-flex
point, by the end of the section we shall indicate how this
hypothesis can be removed. Fix, once and for all, a set of
homogeneous coordinates in $\Bbb P^2$, in such a way that the only
point of $X$ at the infinity line $z=0$ is $(0:1:0)$, this point
will be denote by $\infty$. Thus, $\infty$ is a hyper-flex point:
$X.\{z=0\}=4\infty $.

Denote by $Div^{3,+}(X)$ the set of degree $3$ effective divisors
on $X$, and define the following subset:

$$\align Div^{3,+}_0(X) & := \{D=p_1+p_2+p_3 \in Div^{3,+}(X)\mid h^1(X,\Cal O(D))=0,
\text { and } \infty \ne p_i \} \\
&\cap \{D\in Div^{3,+}(X)\mid p_i\ne \infty\ \text{ and } h^1(X,
\Cal O(p_i+p_j + \infty))=0 \quad \forall i \ne j \}.\endalign
$$

Note that, as the canonical divisor $K_X$ is cut out by the linear
system of lines in $\Bbb P^2$, the above conditions say
geometrically that $\infty$ is not in the support of our divisors
and neither, the three points in the support, nor two points and
$\infty$ are collinear.

Moreover, the condition $h^1(X,\Cal O(D))=0$ allows us to identify
$Div^{3,+}_0(X)$ with a subset of $Pic^3(X)$. Indeed, by
Riemann-Roch any divisor $D\in Div^{3,+}_0(X)$ satisfy $h^0(X,\Cal
O(D))=1$.

The next Theorem explains the basic construction of this section:
we associate to any $D\in Div^{3,+}_0(X)$ the pencil of conics
cutting $X$ out in the divisor $D+\infty$. Our choice of
coordinates gives an explicit form for a basis of this pencil. The
reader must compare this claim with the discussion at the end of
section 2.

\proclaim{Theorem 3.1} There exists bijection between the sets
$Div^{3,+}_0(X)$ and

$$ Z= \{ (A,B,G,H)\in \left( H^0(\Cal O_{\Bbb P^2}(2))\right)^{\oplus 4} \mid F= AG +
BH\},$$ with $A$ and $B$  of the following particular form in
affine coordinates:
$$A= a_{00} + a_{10} x + a_{01}y - x^2,$$

$$B= b_{00} + b_{10} x + b_{01}y - xy,$$
and $G$ satisfying that its coefficient $g_{11}$ corresponding to
the monomial $xy$ equals $1$.
\endproclaim

\demo{Proof} Given $D\in Div^{3,+}_0(X)$, we construct the pair of
conics $A$, $B$ considering the pencil of conics $C$ satisfying
$C.X\ge D+\infty$. In this pencil we fix a basis, namely, that
formed by the only conic $A$ in the system having as tangent line
at $\infty$ the line $\{z=0\}$ and the conic $B$ in the system
defined by the condition $B(1:0:0)=0$.

Sometimes, in order to emphasize the dependence of $A$ and $B$ on
$D$ we write $A_D$, $B_D$ instead of $A$ and $B$.

This conics can be constructed in an explicit way. For instance,
if $p_i\ne p_j$, for $i\ne j$ write $p_i=(x_i,y_i,1)$ (recall that
$p_i\ne \infty$). Consider the matrix:

$$M_D = \left( {\matrix 1 & x_1 & y_1 \\ 1 & x_2 & y_2 \\ 1 & x_3 &
y_3 \endmatrix} \right).$$

This matrix is invertible, because $h^1(X, \Cal O(D))=0$. Thus,
the systems:

$$ M_D \left( {\matrix a_{00} \\ a_{10} \\ a_{01} \endmatrix} \right)=\left( {\matrix x_1^2 \\  x_2^2
 \\  x_3^2 \endmatrix} \right), $$
 and

 $$ M_D \left( {\matrix b_{00} \\ b_{10} \\ b_{01} \endmatrix} \right)=\left( {\matrix x_1y_1 \\
 x_2y_2 \\  x_3y_3 \endmatrix} \right), $$
 have unique solutions. The conics having as coefficients the
 previous solutions are precisely $A$ and $B$.

 If, for example, $D=p_1+2p_2$ we need to modify the system to be
 solved. In this case $A$ is defined by the solutions of:

$$ \left( {\matrix 1 & x_1 & y_1 \\ 1 & x_2 & y_2 \\ 0 & F_y(p_2)
& -F_x(p_2) \endmatrix} \right)\left( {\matrix a_{00} \\ a_{10} \\
a_{01} \endmatrix} \right)=\left( {\matrix x_1^2 \\  x_2^2
 \\  -2x_2F_y(p_2) \endmatrix} \right). $$

 Made this construction, we have the following properties:

 \proclaim{Lemma 3.2} For any $D\in Div^{3,+}_0(X)$, the conics $A$
 and $B$ described above satisfy:

 a) $A$ is irreducible

 b) $a_{01}\ne 0$

 c) $A.B= D+ \infty$

 d) There exists $G$, $H$ $\in H^0(\Bbb P^2, \Cal O_{\Bbb P^2}(2))$ such
 that $F=AG+BH$.

 \endproclaim

 \demo{Proof} a) The equation of $A$ in $\Cal O_{\Bbb P^2,\infty}$
 is given by:

 $$A(x,z)= a_{00}z^2+a_{10}xz-x^2 + a_{01}z. $$

 Assume $A$ is reducible:

$$a_{00}z^2+a_{10}xz-x^2 + a_{01}z= (a+bx+cz)(d+ex+fz),$$
comparing both expressions, we obtain:

$$\align ad & =0 \\ ae+ bd &=0 \\ be & =-1. \endalign $$

If $a=0$ then $bd=0$, but $b=0$ is impossible, because of the
third equation, thus $d=0$. In this way:

$$A=(bx+cz)(ex+fz).$$

As $A.X\ge D$ and the two lines in the factorization contain
$\infty$ the condition $h^1(X,\Cal O(p_i+p_j+\infty))=0$ would be
contradicted.

b) $a_{01}=0$ implies $A$ is reducible.

c) Just by construction.

d) Using Noether's Fundamental Theorem ([5]), it is sufficient to
prove that

$$ F\in <A,B>\Cal O_p,$$
for all $p\in X$. Being $X$ and $A=0$ non-singular, the only non
trivial case to be checked is that of divisors of the form
$D=p_1+2p_2$. In this case,

$$I_{p_i}(F,A)\ge 2= I_{p_i}(A,B),$$
and we are done by a well known criterion ([5], Prop. 1, page
121). \qed
\enddemo

The conics $H$ and $G$ are not uniquely determined by $A$ and $B$,
however, if we consider the Koszul complex:

$$ 0 @>>> \Cal O_{\Bbb P^2} @>(A,B)>> \Cal O_{\Bbb P^2}(2)\oplus \Cal O_{\Bbb
P^2}(2) @>>> \Cal O_{\Bbb P^2}(4) @>>> 0,$$ which is exact and
induces and exact sequence in global sections, we see that all the
possible $G$ and $H$ in the previous relation belong to pencils:

$$G+\lambda B \text{ and } H- \lambda A.$$

Thus, if we fix the unique value of $\lambda$ such that the conic
$G+\lambda B$ have coefficient $g_{11}=1$, then we have defined
uniquely the conics $G$ and $H$ in the relation $F=AG+BH$.

This complete the first part of the proof, namely the assignation
to each $D\in Div_{0}^{+,3}$ of the four conics $(A,B,G,H)$.

Now, we make the converse construction. Assume the conics:

$$A= a_{00} + a_{10} x + a_{01}y - x^2,$$

$$B= b_{00} + b_{10} x + b_{01}y - xy,$$
satisfy the equation

$$F=AG+BH$$
for some conics $G$ and $H$. We shall prove that if $A.B=D+\infty$
then $D \in Div^{3,+}_0(X)$.

Note first that $\infty\notin Supp(D)$, because the intersection
index $I_{\infty}(A,B)=1$. This can be proved by considering the
complete ring ${\hat \Cal O}_{\Bbb P^2,\infty}/<B>$. Looking at
the equation defining $B$ we see that, in this algebra, the class
of $x$ is:

$$\bar x=\frac{b_{00}\bar z^2+ b_{01} \bar z}{1-b_{10} \bar
z}=(1+b_{10}\bar z + b_{10}^2 + ...)(b_{00}\bar z^2+b_{01}\bar
z).$$

Thus, $\bar A \in {\hat \Cal O}_{\Bbb P^2,\infty}/<B>$ can be
expressed like:

$$\bar A(\bar z)= a_{01}\bar z + h(\bar z) \bar z^2.$$

This prove that

$$I_{\infty}(A,B)=\dim_{\Bbb C} {\hat \Cal O}_{\infty}/<A,B>=1.$$

Next, we must show that $h^1(X,\Cal O(D))=0$. If it is not the
case, then there exists a line $L$ such that $X.L\ge D$, but this
implies $A.L\ge D$ and $B.L\ge D$, its follows that $A$ and $B$
have a common factor $L$, contradicting the irreducibility of $F$.

A similar argument implies that $h^1(X,\Cal O(p_i+p_j+\infty))=0$.

This two constructions are clearly one the inverse of the another.
This conclude the proof of the theorem. \qed

\enddemo

The next step in our construction is to prove that the set of
conics previously described has a natural structure of smooth
affine variety, isomorphic with an open set of $Jac X$.

\proclaim{Theorem 3.3} The set $Z$ is a smooth affine variety of
dimension $3$, isomorphic to a Zariski open set of $Jac
X$.\endproclaim

\demo{Proof} {\bf Step 1.} The structure of affine variety.

The structure of affine algebraic set is given exactly by the
relationship:
$$F=AG+BH.$$

More explicitly, $Z$ is by definition a subset $Z\subset \Bbb
A^3\times \Bbb A^3 \times \Bbb A^5\times \Bbb A^6$, with the
coefficients of $(A,B,G,H)$ corresponding to the coordinates
$$((a_{00},a_{10},a_{01}),(b_{00},b_{10},
b_{01}),(g_{00},...)(h_{00},...)).$$

In the set of coordinates corresponding to $G$ there appears only
$5$ entries, because $g_{11}=1$. This set is defined by the
equations

$$F=AG+BH,$$
and being $F$ a fixed curve, this determines  a set of $14$
algebraic equations of degree two in
$(a_{ij},b_{kl},g_{st},h_{uv})$.

This set can be viewed, more accurately, like the image in $\Bbb
A^{17}$ of the map:

$$ \phi: Div^{3,+}_0(X) @>>> Z\subset \Bbb A^{17},$$

$$ D @>>> (A_D,B_D,G_D,H_D).$$

Now, note that two divisors $D,D'\in Div^{3,+}_0(X)$ are linearly
equivalent if and only if they are equal, which follows from the
vanishing of their $h^1$. In other words, $Div^{3,+}_0(X)\subset
Pic^3(X)$. $Div^{3,+}_0(X)$ is in fact an open Zariski set, as can
be proved by noting that in the expression:

$$\align &  Pic^3(X) -Div^{3,+}_0(X)\\ & = \{D=p_1+p_2+p_3 \in Div^{3,+}(X)\mid h^1(X,\Cal O(D))=1,
\text { or } \infty = p_i \} \\
&\cup \{D\in Div^{3,+}(X)\mid h^1(X, \Cal O(p_i+p_j + \infty))=1
\text { for some } i \ne j \text{ or } p_i = \infty\},\endalign $$
the first set is just the image of $W_2\subset Pic^2(X)$ under the
map $+\infty: Pic^2(X)@>>> Pic^3(X)$, and the second can be
obtained as the epimorphic  image of the morphism:

$$\pi: X\times X @>>> Pic^3(X),$$
defined as follows: given $(p,q)$, be $L$  the line determined by
$q$ and $\infty$, and set $p_2+p_3= X.L-q-\infty$. Then, define
$\pi (p,q)$ as the class of $p+p_2+p_3$.

Moreover, $\phi$ is easily expressed as a function of the
symmetric functions on $x_i$ and $y_i$. In this way, $\phi$ is an
injective algebraic map onto $Z$. In particular, $Z$ is
irreducible.

{\bf Step 2.} Dimension of $Z$.

Let us to compute the Zariski dimension of $Z$ at any point and to
prove it is equal to $3$. Note that $Z$ is defined by $14$
equations in $\Bbb A^{17}$.

Given $(A,B,G,H)\in Z$, consider the tangent space:

$$\align T_{(A,B,G,H)} Z=& \{ (\dot A, \dot B,\dot G,\dot H)\mid deg \dot A, \dot B, \dot G, \dot H \le 2 \text { and} \\
& F=(A+\epsilon \dot A)(G+\epsilon \dot G)+ (B+\epsilon \dot
B)(H+\epsilon \dot H) \mod \epsilon^2\} \endalign$$

Thus, we must have:

$$A\dot G+\dot A G+ B\dot H + \dot B H=0. \tag 2.1$$

Hence, our problem reduces to prove that any quartic in $\Bbb P^2$
containing the point $(0:1:0)$ can be expressed in the form 2.1,
for some $\dot A, \dot B, \dot G, \dot H$.

Let $V=\{A\dot G+\dot A G+ B\dot H + \dot B H\}$ and write

$$V= V_1 \oplus V_2,$$
with $V_1= \{ \dot A =\dot B = \dot G =0 \}$ and $V_2=\{\dot
H=0\}$.

The space $V_1$ equals the space of quartics

$$\{ Q\mid B\text{ divides } Q\}\simeq H^0(\Bbb P^2, \Cal O_{\Bbb P^2}(2)),$$
hence, $\dim V_1=6$.

Moreover, we have that inside $V_2$ the subspaces $V_2^0\subset
V_2$, defined by $V_2^0=\{ \dot A=\dot B=0\}$ also have dimension
$6$. Thus, we only need to prove that, inside $V_2$, there exists
a complementary 2-dimensional subspace to $V_2^0$.

It is easily provided by the space generated by $H$ and $G$, that
will be complementary to $V_2^0$, unless $sH+tG=A$ for some
constant $s$ and $t$.

If this were the case, we argue in a similar fashion using now the
decomposition:

$$V= \{ \dot A= \dot B=\dot H=0\} \oplus \{ \dot G=0 \},$$
and use that $sH+tG=A$ and $aG+bH=B$ implies $F$ is reducible.
\qed
\enddemo

In this way, we have explicitly constructed an affine variety
structure for an open set of $Pic^3(X)$.

{\bf Remarks} 1. Being this construction only a concrete
realization of the general discussion in section 2, it seems
obvious that the hypothesis on the existence of an hyper-flex is
superfluous. In fact, a similar construction can be carried out if
we eliminate this hypothesis, only small modifications are needed
in the definition of $Div^{+,3}_0(X)$, fixing a point having as
tangent $\{z=0\}$ and avoiding divisors containing in its support
any point on $z=0$.

2. It can be proved, in the spirit of [7], that a finite number of
 translations of $Z$ give an affine covering of $Pic^3(X)$,
 recovering the structure of $Pic^3(X)$ as a complete variety. The
 divisors needed for this have support on the ramifications points
 of the projection of $X$ with center at $\infty$ on a suitable
 projective line. The details can be consulted in [8].

 \head 4. Explicit law group on $JX$. \endhead

 Let $X$ be a canonical curve of genus $3$, assume for simplicity that $X$ has
 a hyper-flex point $\infty \in X$.

 In this section we describe explicitly the law
 group on $JX$. The description includes an explicit algorithm for
 computing the inverse of an element in $JX$ and the
 sum of two elements. The computation of the inverse is
 particulary transparent when restricted to the open set $Z$
 defined in the previous section. To start with, consider the
 morphism:

 $$ \pi_{\infty}: Pic^3(X)  @>>> JX $$
 $$D  @>>> D-3\infty . $$

 It is an elementary consequence of Riemann-Roch Theorem (and also
 a part of Jacobi Inversion Theorem) that this map is surjective.
 Moreover, as previously noted, it is injective in
 $Div^{3,+}_0(X)$. By a slight abuse of notation we denote
 indistinctly by $Z$ to both, the image in $JX$ of $Div^{3,+}_0(X)$
 and to the affine variety constructed in section 2 (they are, in
 fact, isomorphic).

 \proclaim{Theorem 4.1} The open set $Z$ is invariant under
 multiplication by $(-1)$ in $JX$. In the affine coordinates for
 $Z$ this map is given by:

 $$ (A,B,G,H) @>>> (A,H,G,B).$$

\endproclaim

\demo{Proof} Let $D\in Div^{3,+}_0(X)$, denote by $D^-$ the
divisor such that $D+D^--6\infty \equiv 0$. This divisor can be
constructed as follows: let $A$ be the conic considered before,
then $X.A=D+ D'+2\infty$ with $deg D'=3$. Let
$L_{\infty}=\{z=0\}$, being $\infty$ a hyper-flex we have:

$$X.A\equiv X.2L_{\infty}\equiv 8\infty,$$
thus, $D'=D^-$.

Let us prove that $D^-\in Div^{3,+}_0(X)$. A computation similar
to the local argument invoked in the proof of Theorem 2.1 shows
that $I_{\infty}(X,A)=2$, thus $\infty \notin Sup D'$. Moreover,
if $h^1(X,\Cal O_X(D'))=1$ then $A$ must be reducible and the same
is valid if $h^1(X,\Cal O_X(q_i+q_j+\infty))=1$ with $q_i$, $q_j$
points in the support of $D'$. This prove that $Z$ is invariant
under multiplication by $-1$.

Now, the conic $A_D$ coincides with $A_{D^-}$, as previously
noted. Finally, as

$$F=AG_D+B_DH_D=AG_{D^-}+ B_{D^-}H_{D^-},$$
we have
$$\align A.F=D+D^-+2\infty  & = A.B_{D^-}+A.H_{D^-} \\
& = D^-+\infty + A. H_{D^-} \endalign $$ and thus

$$A.B_D=A.H_{D^-}.$$
By the uniqueness of the decomposition $F=AG+BH$, once fixed $A$
and $B$, we have $B_{D^-}=H_D$.

In this way the set of conics associated to $D^-$ is
$(A_D,H_D,G_D,B_D)$. \qed
\enddemo

We remark that a conic $A$ satisfying $A.X\ge D+ 2\infty$ exists
independently of being $D\in Div^{3,+}_0(X)$ or not. Thus, an
explicit and simple algorithm for computing the inverse of a
divisor $D-3\infty$ is to take the class of the divisor
$A.X-D-5\infty$.

The proof also shows that the conic $H_D$ has coefficients
$h_{20}=h_{02}=0$. We obtain as a by-product:

\proclaim{Theorem 4.2} The open set of the Kummer variety $K(X)$
given by the quotient $Z/{\{\pm 1\}}$ is determined by the
equations:

$$ F=A.G+ Q,$$
with $Q$ a quartic expressable like a product of two conics
passing through $\infty$ and $(1:0:0).$ \endproclaim

The condition on the reducibility of the quartic $Q$ can be
written down explicitly in terms of the usual procedure of
composing the Segre map with a linear projection. More concretely,
a quadric $Q$ will be the product of such two conics if when
expressed as

$$\sum_{ij} c_{ij} x^i y^j,$$
its coefficients satisfy:

$$c_{30}=c_{03}=c_{40}=c_{13}=c_{31}=c_{04}=0,$$

$$\align c_{22} =1 \quad c_{00} =w_{00}, \quad c_{10} =w_{12}+w_{21} \\ c_{01} =
w_{20}+w_{02} \quad c_{11} =w_{03}w_{30}+w_{12}w_{21}\quad c_{20}
&= w_{21} \\ c_{02}  =w_{12} \quad c_{21} = w_{13}+w_{31} \quad
c_{12}  = w_{23} + w_{32} \quad c_{22} &=w_{33}, \endalign $$ and
$w_{ij}$ satisfying the Segre relations $rank (w_{ij})\le 2$.

\vskip0.2cm

The addition law on $JX$, considered like the image of
$\pi_{\infty}$, can be also described by a simple algorithm:

Let $D-3\infty, \quad D'-3\infty \in JX$. Consider a cubic $C$
such that $C.X= D+ D'+ E^{-} +3\infty$.

Let $A_{E^-}$ be a conic with $A_{E^-}.X=E^- + E+2\infty$.

Now, $D+D'+E^-\equiv 9\infty$ and

$$ E^-+E\equiv 6\infty.$$
It follows that $D+D'-6\infty \equiv E+ 3\infty,$ that is

$$ D+D'\equiv E-3\infty.$$

Of course, $Z$ can not be invariant under the addition law (it can
not be a subgroup of $JX$ for several reasons), so we can not
expect to find a description of the sum in terms of the affine
coordinates for $Z$.

\head 5. The description via Grassmanians \endhead

In this section we outline a related construction in terms of
Grassmanians. This construction is a bite more intrinsic, as it
does not involve the particular choice of $A$ and $B$ made before.

Let $V= H^0(\Bbb P^2, \Cal O_{\Bbb P^2}(2)\otimes \frak
m_{\infty})$. The construction in section 3 assigns to any $D \in
Div^{3,+}_0(X)$ a line in $\Bbb PV$, i.e., the pencil $\lambda A+
\mu B$.

We have then an injective morphism:

$$ Z\cong Div^{3,+}_0(X) @>>> \Bbb G (1,\Bbb P V),$$

$$ (A,B,G,H) @>>> \lambda A+ \mu B .$$

In terms of the Pl\" uckers embedding of the Grassmanian this
morphism can be described explicitly. If we call

$$A= a_0 + a_1 x + a_2y - x^2,$$

$$B= b_0 + b_1 x + b_2y - xy,$$
then the Pl\"ucker embedding restricted to $Z$ is the map sending
$(A,B)$ into
$$((a_1b_0-b_1a_0, a_2b_0-b_2a_0,-a_0, a_1b_2-a_2b_1, -a_1,-a_2,
-b_0,-b_1,-b_2,1)).$$

This together with Pl\" ucker relations give another set of
equations for $Z$. Note that it is obvious from the expression for
the map that this is not only injective, but in fact an embedding
when restricted to $Z$.

It is worthy to be noted that this map can be extended to any
divisor $D\in Pic^3(X)$ such that $h^1(X, \Cal O_X(D))=0$, as
under this condition a pencil of conics cutting  $X$ out on $D+
\infty$ is well defined. In other words, we can think of a
morphism:

$$ JX-\Theta @>>> \Bbb G (1, \Bbb P V).$$
This morphism, however, will be injective only on $Z$.


\Refs

\ref \no 1 \paper An Explicit Algebraic Representation of the Abel
Map \by G. W.  Arderson\jour  Int. Math. Research Notices \vol
11\page 495--205 \yr 1997
\endref

\ref \no 2 \paper Jacobians of spectral curves and completetly
integrable Hamiltonian system  \by A. Beauville\jour Acta
Mathematica \vol 164\page 211--235 \yr 1990
\endref

\ref \no 3 \by F. Catanese \paper Homological algebra and
algebraic surfaces \jour Algebraic geometry Santa Cruz 1995, Proc.
Sympos. Pure Math. Amer. Math. Soc., Providence, RI.\yr 1997 \vol
62 \page 3--56 \endref

\ref \no 4 \paper On the Jacobian varieties of Picard curves:
addition law and algebraic structure\by J. Estrada-Sarlabous, E.
Reinaldo Barreiro, and A. Pi\~neiro Barcel\'o \jour M. Nachrichten
 \vol 208\page 149--166 \yr 1999
\endref

\ref \no 5 \by W. Fulton \book Algebraic Curves \publ
Benjamin-Cunning, \yr 1969
\endref

\ref \no 6 \by M. Green \paper Koszul cohomology and geometry
\jour ICTP College "Lectures on Riemann Surfaces" World Scientific
Press \yr 1989  \page 177--200
\endref

 \ref \no 7 \by D. Mumford \book Tata
Lectures on Theta II \publ Modern Birkh\"auser Classics \yr 2007
\endref

\ref \no  8\book Construcci\'on algebraica de la Jacobiana de una
curva de g\'enero $3$ \publ Ph. D. Thesis, CIMAT \by J.
Romero-Valencia \yr 2008 \endref

\ref \no 9 \paper Some remarks on the Jacobian varieties of a
Picard Curve \by A. G. Zamora\jour Bolet\'{\i}n SMM \vol 3 \page
221--234 \yr 1997
\endref

\endRefs
\enddocument